\documentclass[12pt,a4paper]{amsart}
\usepackage[english]{babel}
\usepackage[utf8]{inputenc}
\usepackage{fancyhdr}
\usepackage{xcolor}

\usepackage{amssymb}
\usepackage{amsmath} 
\usepackage{amsthm} 
\usepackage{enumitem}
\usepackage{bbm}

\newtheorem{thm}{Theorem}
\newtheorem{prop}{Proposition}
\newtheorem{lemma}{Lemma}
\newtheorem*{thm*}{Theorem}

\theoremstyle{proposition}
\newtheorem*{prop*}{Proposition}
\theoremstyle{proposition}
\newtheorem*{rem*}{Remark}

\numberwithin{lemma}{section}
\numberwithin{equation}{section}
\numberwithin{prop}{section}
\numberwithin{thm}{section}

\title{A weighted central limit theorem for $\log|\zeta(1/2+it)|$}
\author{Alessandro Fazzari}
\address{Universit\`a  di Genova,  Dipartimento di Matematica. 
 Via Dodecaneso 35, 16146 Genova, Italy}
\email{fazzari@dima.unige.it}

\DeclareMathOperator{\meas}{meas}

\begin{document}
\maketitle
\begin{abstract} Under the Riemann Hypothesis, we show that as $t$ varies in $T\leq t \leq 2T$, the distribution of $\log|\zeta(1/2+it)|$ with respect to the measure $|\zeta(1/2+it)|^2dt$ is approximately normal with mean $\log\log T$ and variance $\frac{1}{2}\log\log T$.\end{abstract}

\section{Introduction and statement of the main results}

In the understanding of the value distribution of the Riemann zeta function on the critical line, the milestone is due to Selberg, who proved a central limit theorem for $\log |\zeta(1/2+it)|$, showing that 
\begin{equation}\label{CLT} \frac{1}{T} \meas \Bigg \{t\in[T,2T]:\frac{ \log |\zeta (\frac{1}{2}+it)|}{\sqrt{\frac{1}{2}\log\log T}} \geq V  \Bigg \}\sim \int_V^{\infty} e^{-\frac{x^2}{2}}\frac{dx}{\sqrt{2\pi}} \end{equation}
for any fixed $V$, as $T$ goes to infinity (see \cite{Sel} for the first results in this direction and Tsang's thesis \cite{Tsang} for an actual formulation of the central limit theorem). Analogous statements hold also in more generality, for example in the case of the imaginary part of $\log\zeta(1/2+it)$ or other $L-$functions (see e.g. \cite{Sel2} or \cite{B-H}). In 2015 Radziwill and Soundararajan \cite{RS} gave a new and simple proof of Selberg's central limit theorem. \\
In this context one may ask about the uniformity in $V$ of (\ref{CLT}), investigating the large values of zeta. Classically it was known that (\ref{CLT}) holds for $V=V(T)\ll (\log\log\log T)^{1/2-\varepsilon}$, $\varepsilon>0$ (see \cite{Tsang}). More recently Radziwill \cite{Radz} introduced a new method that extended (\ref{CLT}) to the large deviation range $V \ll (\log\log T)^{1/10-\varepsilon}$. Furthermore he conjectured that the largest range of uniformity for (\ref{CLT}) is $V=o(\sqrt{\log\log T})$. Recently Inoue \cite{Shota} proved a partial result towards this conjecture.\\ 
Moreover Soundararajan \cite{Sound} obtained upper bounds for the measure of the set of large values, in the case $V\gg \sqrt{\log\log T}$ (and an unconditional lower bound in \cite{Sound2}). 
The speculation that gives rise to Soundararajan's work is that an upper bound like 
\begin{equation}\label{conj} \frac{1}{T}\meas \Bigg \{t\in[T,2T]: \frac{\log|\zeta (1/2+it)|}{\sqrt{\frac{1}{2}\log\log T}} \geq V  \Bigg \} \ll \frac{1}{V}\exp\Big(-\frac{V^2}{2}\;\Big)\end{equation}
also holds for some $V\gg \sqrt{\log\log T}$. Even though Soundararajan does not prove such a precise upper bound, he gets a quasi optimal one under the Riemann Hypothesis, which is enough to derive conditional quasi optimal upper bounds for the moments of the Riemann zeta function. Then the problem of proving (\ref{conj}) is still open. For example in the case $V=\sqrt{2\log\log T}$, Soundararajan only proved that the left hand side of (\ref{conj}) is essentially $ \ll(\log T)^{-1+o(1)} $ while the conjectural sharp upper bound should be $ \ll (\log T \sqrt{\log\log T})^{-1} $ (see Harper \cite{Har1}, \cite{Har2} and \cite{Har3} for further discussions). \\

In this setting, with the aim of studying the large values of the Riemann zeta function, one can investigate the distribution of $\log |\zeta(1/2+it)|$ with respect to a different measure. For instance one can ``tilt'' the measure and study the distribution of $\log|\zeta (1/2+it)|$ with respect to the weighted measure
$$ |\zeta|^2dt:=|\zeta(1/2+it)|^2dt.$$
This change of measure means that in integrals which represent probabilities (or moments) we are giving more importance to the contribution of those $t$ such that $|\zeta(1/2+it)|$ is large. For this reason, understanding the distribution of $\log|\zeta(1/2+it)|$ with respect to the weighted measure $|\zeta|^2dt$ might be of help in the understanding of the large values of $\zeta$. 
\begin{thm}\label{thm1}
Under the Riemann Hypothesis, as $t$ varies in $T\leq t\leq 2T$, the distribution of $\log|\zeta(1/2+it)|$ is asymptotically Gaussian with mean $\log\log T$ and variance $\frac{1}{2}\log\log T$, with respect to the
weighted measure $|\zeta|^2dt$.
\end{thm}

We note that this result is a manifestation of Girsanov's theorem from probability theory, which describes how a stochastic process changes under certain changes of measure. In the specific case we are interested in, Girsanov's theorem reduces to the simple fact that if we take $X$ a Gaussian random variable with mean $0$ and variance $\sigma^2$ and we tilt $X$ against $e^{yX}$ with $y\in\mathbb R$, then the resulting random variable is again Gaussian, with the same variance but mean $y\sigma^2$ (it can be proved just by completing the square). Theorem 1 shows the same behavior for $\log|\zeta(1/2+it)|$ in the case $y=2$, reinforcing our expectation that $\log|\zeta(1/2+it)|$ behaves like a Gaussian in many respects (other interesting computations involving the Riemann zeta function inspired by Girsanov's theorem can be found in Harper's work \cite{Har4}, Section 3).\\

We now describe the general strategy to prove Theorem \ref{thm1}. Even though the Euler product formula only holds in the half-plane of convergence, for many purposes the Riemann zeta function behaves like an Euler product also on the critical line (see Principle 1.3 in \cite{Har3}), thus $\log|\zeta(1/2+it)|$ behaves like a Dirichlet polynomial. Roughly speaking we know that for a suitable $x=x(T)$ we have
\begin{equation}\label{roughly}\log|\zeta(1/2+it)|\approx\Re\sum_{p\leq x}\frac{1}{p^{1/2+it}}+(\text{contribution from zeros})\end{equation}
 (see \cite{Har1}, \cite{Har2} for further and more precise details) and in several applications the contribution from the zeros can be controlled (two important examples are \cite{Sound} and \cite{Har1}).
This approximation also holds in our setting, as shown by the following proposition:
\begin{prop}\label{prop1}
Let $T$ be a large parameter. Denote $P(t)= \sum_{p\leq x}p^{-1/2-it}$, where $x=T^{\varepsilon/k}$, $\varepsilon:=(\log\log\log T)^{-1}$, $k$ a positive integer. Under the Riemann Hypothesis, there exists a constant $C>0$ such that we have uniformly in $k$:
$$ \frac{1}{T\log T}\int_T^{2T}\big | \log |\zeta (1/2+it)|-\Re P(t)\big |^{2k}|\zeta|^2dt =O\left ( (Ck)^{4k}(\log\log\log T)^{2k+1/2} \right ). $$
\end{prop}
We remark that this is the only point where we rely on the assumption of the Riemann Hypothesis. In fact in order to estimate the contribution of the zeros that appears in (\ref{roughly}), we need to bound the sum over the non-trivial zeros $ \sum_{0<\rho \leq T}|\zeta(\rho+i\alpha)|^2 $ with $|\alpha|\leq 1$ a real parameter, which is known to be $\ll T(\log T)^2$ only conditionally on the Riemann Hypothesis (see \cite{Gonek}).  \\

Thanks to Proposition \ref{prop1}, at this point it suffices to show that the distribution of $\Re P(t)$ is approximately Gaussian with respect to the measure $|\zeta|^2dt$. This is achieved by the method of moments using the following result.
\begin{prop}\label{prop2}
Let $P(t)=\sum_{p\leq x}p^{-1/2-it}$, $x:=T^{\varepsilon/k}$, $\varepsilon:=(\log\log\log T)^{-1}$. Denote $\mathcal L=\sum_{p\leq x}\frac{1}{p}$. Then, for every fixed $k$ integer 
$$\frac{1}{T\log T}\int_T^{2T}(\Re P(t)-\mathcal L)^k|\zeta|^2dt=\begin{cases} \big (\frac{\mathcal L}{2}\big )^{k/2}(k-1)!!+O_k\big(\mathcal L^{(k-1)/2}\big) \hspace{0.75cm} \text{if } k \text{ is even} \\ O_k\big(\mathcal L^{(k-1)/2}\big) \hspace{3.9cm}\text{if } k \text{ is odd}.\end{cases}$$
\end{prop}
Note that by definition of $x$ we know that $\log\log x=\log\log T-\log k+\log\varepsilon$, then for a fixed $k$ we have $\log\log x= \log\log T+O(\log_4T)$, where $\log_4$ denotes the fourth iterated natural logarithm. Hence by Mertens' theorem $\mathcal L=\log\log x+O(1)=\log\log T+O(\log_4 T)$. As a consequence the right hand side in Proposition \ref{prop2} matches with the moments of a normal of mean $\mathcal L\sim\log\log T$ and variance $\frac{\mathcal L}{2}\sim\frac{1}{2}\log\log T$. Putting together the two propositions one has that the moments of $\log|\zeta(1/2+it)|$ with respect to the measure $|\zeta|^2dt$ are asymptotic to the moments of a Gaussian random variable of mean $\log\log T$ and variance $\frac{1}{2}\log\log T$, then the theorem is proved. \\

Finally, the author would like to remark that a similar analysis can be performed also considering the tilted measure $|\zeta|^4dt$ and, under suitable assumptions, the general case $|\zeta|^{2k}dt$ with $k\in\mathbb N$. This problem will be addressed in \cite{2.}.


\section{Proof of Proposition \ref{prop1}}
First of all, we state an important tool which allows us to compute the moments of a sufficiently short Dirichlet polynomial with respect to $|\zeta|^2dt$, which can be established by the method in \cite{BCHB} and \cite{BCR}. 
\begin{lemma}\label{B-C-HB}
Let $A(s)=\sum_{n\leq T^\theta}a(n)n^{-s}$ and $B(s)=\sum_{m\leq T^\sigma}b(m)m^{-s}$ be Dirichlet polynomials with $a(n)\ll n^\varepsilon$, $b(m)\ll m^\varepsilon$ for every $\varepsilon >0$ and $\theta+2\sigma< 1$. Then, denoting $c:=2\gamma +\log 4 -\log 2\pi-1$, we have:
$$ \int_T^{2T}A(1/2+it)|B(1/2+it)|^2|\zeta|^2dt = T\sum_{m,n}\frac{(a*b)(n)\overline{b(m)}}{[n,m]}\bigg (\log\bigg( \frac{T(n,m)^2}{nm}\bigg)+c \bigg ) + o(T) $$
where $(a*b)(n)=\sum_{d_1d_2=n}a(d_1)b(d_2)$ denotes the standard convolution.
\end{lemma}
Our proof of Proposition \ref{prop1} is a modification of Theorem 5.1 in \cite{Tsang}. We recall that $P(t)=\sum_{p\leq x}p^{-1/2-it}$ and $x=T^{\varepsilon/k}$ with $\varepsilon=(\log\log\log T)^{-1}$. Following Tsang's strategy, whose notations become easier under the Riemann Hypothesis, we have (see \cite{Tsang}, equation (5.15)):
\begin{equation}\label{approx1} \log \zeta(1/2+it)-P(t) = S_1 +S_2 +S_3 + O(R) - L(t)  \end{equation}
with 
$$ S_1:= \sum_{p\leq x}\big ( p^{-1/2-4/\log x}-p^{-1/2}\big )p^{-it},\hspace{1cm} 
S_2:=\sum_{\substack{p^r\leq x \\ r\geq 2}}\frac{p^{-r(1/2+4/\log x+it)}}{r}$$
$$ S_3:=\sum_{x<p\leq x^3}\frac{\Lambda(n)}{\log n}n^{-1/2-4/\log x-it}, \hspace{1cm}
 R:= \frac{5}{\log x}\bigg ( \bigg | \sum_{n\leq x^3}\frac{\Lambda (n)}{n^{1/2+4/\log x+it}} \bigg |+\log T \bigg ) $$
$$ L(t):=\sum_\rho \int_{1/2}^{1/2+4/\log x}\left ( \frac{1}{2}+\frac{4}{\log x}-u \right )\frac{1}{u+it-\rho}\frac{1}{\frac{1}{2}+\frac{4}{\log x}-\rho}du,$$
where the sum in the definition of $L(t)$ is over all the non-trivial zeros of $\zeta$. Hence
\begin{equation}\label{approx2} \log |\zeta(1/2+it)|-\Re (P(t)) = \Re (S_1) +\Re (S_2) +\Re (S_3) + O(R) -\Re ( L(t))  \end{equation}
so what remains to do is studying the $2k-$moments of all these objects with respect to the weighted measure $|\zeta|^2dt$; to this aim we rely on Lemma \ref{B-C-HB}.\\

Let's start with the first one:
\begin{equation}\begin{split}\notag
\frac{1}{T\log T}\int_T^{2T} |S_1|^{2k} |\zeta|^2dt&=\frac{1}{T\log T} \int_T^{2T} \bigg | \bigg ( \sum_{p\leq x}\frac{p^{-4/\log x}-1}{p^{1/2+it}} \bigg )^k \bigg |^2|\zeta|^2dt 
\\&=  \frac{1}{T\log T}\int_T^{2T} \bigg | \sum_{n\leq x^k}\frac{1}{n^{1/2+it}}\sum_{\substack{p_1\cdots p_k=n \\ p_i \leq x \: \forall i}}\prod_{i=1}^k\big ( p_i^{-4/\log x}-1 \big ) \bigg |^2|\zeta|^2dt 
\\& \ll \sum_{m,n\leq x^k}\frac{(m,n)}{mn}\sum_{\substack{p_1\cdots p_k=n \\ q_1\cdots q_k=m \\ p_i,q_i\leq x}} \prod_{i=1}^k|p_i^{-4/\log x}-1| |q_i^{-4/\log x}-1|
\\& =  \sum_{\substack{p_1,\dots ,p_k\leq x \\ q_1,\dots q_k\leq x}}\frac{(p_1\cdots p_k,q_1\cdots q_k)}{p_1\cdots p_k q_1\cdots q_k}\prod_{i=1}^k|p_i^{-4/\log x}-1||q_i^{-4/\log x}-1|.
\end{split}\end{equation}
To make the GCD on the numerator explicit, we rewrite the primes $p_1,\dots,p_k$ highlighting the multiplicity of these primes:
$$ \{ p_1,\dots,p_k \}=\{p_1',\dots, p_l' \}  $$
where the $p_i'$'s are distinct and we denote $c_i\geq 1$ the multiplicity of $p_i'$ in this set, so $c_1+\dots+c_l=k$. Now we do the same for the $q_i$'s and we put in evidence if any $q_i$ already appears among the $p_i'$'s:
$$ \{ q_1,\dots,q_k \} = \{ p_1',\dots ,p_l' \} \cup \{q_1',\dots ,q_m'\} $$ 
where the $p_i'$'s and $q_j'$'s are all distinct and we denote $e_i\geq 0$ and $d_i\geq 1$ the multiplicities of $p_i'$ and $q_i'$ respectively. Then we have $e_1+\dots +e_l+d_1+\dots +d_m=k$. In the following we drop the symbol $'$, just denoting the new primes with $p_i,q_i$. 
With these notations, the previous sum is 
$$ \ll (k!)^2 \sum_{\substack{l\leq k \\ m\leq k}}  \sum_{\substack{c_1+\dots c_l=k \\ e_1+\dots e_l+d_1+ \dots +d_m=k \\ c_i\geq 1, \; d_i \geq 1, \; e_i \geq 0}} \prod_{i=1}^l \Bigg ( \sum_{p_i} \frac{|p_i^{-4/\log x}-1|^{c_i+e_i}}{p_i^{\max (c_i,e_i)}} \Bigg ) \prod_{i=1}^m \Bigg (\sum_{q_i}\frac{|q_i^{-4/\log x}-1|^{d_i}}{q_i^{d_i}} \Bigg )$$
and if we ignore the equation for $c_i,e_i,d_i$ we get 
\begin{equation}\label{vhs}  \ll (k!)^2 \sum_{\substack{l\leq k \\m\leq k}}\prod_{i=1}^l\Bigg ( \sum_{\substack{c_i\geq 1\\ e_i\geq 0}}\sum_{p_i\leq x}\frac{|p_i^{-4/\log x}-1|^{c_i+e_i}}{p_i^{\max(c_i,e_i)}} \Bigg ) \prod_{i=1}^m \Bigg (\sum_{d_i\geq 1}\sum_{q_i\leq x}\frac{|q_i^{-4/\log x}-1|^{d_i}}{q_i^{d_i}}\Bigg ).  \end{equation}
Now we remark that only in the case $c_i=1$ and $e_i\leq1$ the sum over $p_i$ in the first parentheses gives an unbounded contribution. Indeed the remaining cases give 
\begin{equation}\begin{split}\notag 
\sum_{\substack{c_i\geq 1,\;e_i\geq 0: \\ \max(c_i,e_i)\geq 2}}& \sum_{p_i\leq x}\frac{|p_i^{-4/\log x}-1|^{c_i+e_i}}{p_i^{\max(c_i,e_i)}}
\ll  \sum_{\substack{c_i\geq 1,\;e_i\geq 0: \\ \max(c_i,e_i)\geq 2}} \sum_{p_i\leq x}\frac{1}{p_i^{\max(c_i,e_i)-3/2}}\frac{1}{p_i^{3/2}}
\\&\ll \sum_{\substack{c_i\geq 1,\;e_i\geq 0: \\ \max(c_i,e_i)\geq 2}}\frac{1}{2^{\max(c_i,e_i)-3/2}} \sum_{p_i\leq x}\frac{1}{p_i^{3/2}}\ll \sum_{\substack{c_i\geq 1,\;e_i\geq 0: \\ \max(c_i,e_i)\geq 2}}\frac{1}{2^{\max(c_i,e_i)}}\ll \sum_{\substack{c_i\geq 0\\e_i\geq 0}}\frac{1}{2^{(c_i+e_i)/2}}\ll 1.
\end{split}\end{equation}
We treat the second parentheses analogously, so that we get a bound for (\ref{vhs}), which is
\begin{equation}\begin{split}\notag
&\ll (k!)^2 \sum_{\substack{l,m\leq k}}\bigg ( \sum_{p\leq x}\frac{|p^{-4/\log x}-1|}{p}+\sum_{p\leq x}\frac{|p^{-4/\log x}-1|^2}{p}+ O(1)\bigg )^{l+m}.
\end{split}\end{equation}
In order to bound the first sum we use that $e^{-z}=1+O(z)$ for $z\ll 1$, so if $p\leq x^\delta$ (with $\delta\ll 1$) we have  $$ \left | p^{-4/\log x}-1\right | = \left |e^{-\frac{4\log p}{\log x}}-1 \right | \ll \frac{4\log p}{\log x}\ll \delta$$ and as a consequence 
$$ \sum_{p\leq x^\delta}\frac{\left | p^{-4/\log x}-1\right |}{p}\ll \delta \sum_{p\leq x^\delta}\frac{1}{p}\ll \delta  \big (|\log \delta|+\log\log x\big ). $$
On the other hand, if $x^\delta <p \leq x$ than trivially $|p^{-4/\log x}-1|\leq 2$ hence
$$ \sum_{x^\delta<p\leq x}\frac{\left |p^{-4/\log x}-1\right |}{p}\ll \sum_{x^\delta<p\leq x}\frac{1}{p}\ll |\log \delta|+O(1). $$
Therefore  
$$ \sum_{p\leq x} \frac{\left | p^{-4/\log x}-1\right |}{p} \ll \delta  \big ( |\log \delta|+\log\log x \big ) + |\log \delta|+O(1) \ll \log\log\log T$$ 
by selecting $\delta = \frac{\log\log\log x}{\log\log x}$. The second sum is $\ll \log\log\log T$ too, being $|A-1|^2\leq |A-1|$ for $0<A<1$.
Putting all together, the sum we are considering is:
$$  \ll (k!)^2 \sum_{\substack{l,m\leq k}}(3\log\log\log T)^{l+m} \ll (k!)^2 C^{2k} (\log\log\log T)^{2k}$$
with $C$ a sufficiently large positive constant. In conclusion, uniformly in $k$ we have  $$ \frac{1}{T\log T}\int_T^{2T} |S_1|^{2k} |\zeta|^2dt  \ll (Ck)^{2k} (\log\log\log T)^{2k}.$$

Now we focus on $S_2$. Using again Lemma \ref{B-C-HB} we have:
\begin{equation}\begin{split}\notag
\frac{1}{T\log T}&\int_T^{2T} |S_2|^{2k} |\zeta|^2dt 
\ll  \sum_{\substack{p_1,\dots,p_k \leq x \\ q_1,\dots,q_k\leq x \\ r_1,\dots,r_k\geq 2 \\ s_1,\dots,s_k\geq 2}}\frac{(p_1^{r_1}\cdots p_k^{r_k}, q_1^{s_1}\cdots q_k^{s_k})}{p_1^{r_1}\cdots p_k^{r_k} q_1^{s_1}\cdots q_k^{s_k}}.
\end{split}\end{equation}
We use the same decomposition of $\{p_1,\dots,p_k\}$ and $\{q_1,\dots,q_k\}$ as before getting 
\begin{equation}\begin{split}\notag
&\ll(k!)^2\sum_{m,l\leq k}\sum_{\substack{p_1,\dots,p_l\leq x \\q_1,\dots,q_m\leq x}}\sum_{\substack{a_1,\dots,a_l\geq 2 \\ b_1,\dots,b_l\geq 0 \\ f_1,\dots,f_m\geq 2}}\frac{1}{p_1^{\max(a_1,b_1)}\cdots p_l^{\max(a_l,b_l)}q_1^{f_1}\cdots q_m^{f_m}} \\ 
& =(k!)^2 \sum_{m,l\leq k}\prod_{i=1}^l \Bigg ( \sum_{\substack{p_i \leq x \\ a_i\geq 2 \\ b_i \geq 0}} \frac{1}{p_i^{\max(a_i,b_i)}}\Bigg ) \prod_{i=1}^m \Bigg ( \sum_{\substack{q_i\leq x \\ f_i \geq 2}}\frac{1}{q_i^{f_i}} \Bigg )
\ll (k!)^2\sum_{m,l\leq k}C^{m+l} \ll (Ck)^{2k}.
\end{split}\end{equation}

Let us investigate $S_3$ using the same approach. We have
\begin{equation}\begin{split}\label{dvd}
\frac{1}{T\log T}\int_T^{2T} &|S_3|^{2k}|\zeta|^2dt\ll \sum_{\substack{x<p_1^{r_1},\dots,p_k^{r_k}\leq x^3 \\ x<q_1^{s_1},\dots,q_k^{s_k}\leq x^3 }}\frac{(p_1^{r_1}\cdots p_k^{r_k},q_1^{s_1}\cdots q_k^{s_k})}{p_1^{r_1}\cdots p_k^{r_k}q_1^{s_1}\cdots q_k^{s_k}}.
\end{split}\end{equation}
We begin studying the case when all the exponents $r_i,s_i$ are equal to $1$. We can implement the same technique as before, getting:
\begin{equation}\begin{split}\notag
\ll(k!)^2\sum_{m,l\leq k}&\sum_{\substack{c_1,\dots,c_l\geq 1 \\ e_1,\dots,e_l\geq 0 \\ d_1,\dots,d_m \geq 1}} \sum_{\substack{x<p_1,\dots,p_l\leq x^3 \\ x<q_1,\dots,q_m\leq x^3}}\frac{1}{p_1^{\max(c_1,e_1)}\cdots p_l^{\max(c_1,e_l)}q_1^{d_1}\cdots q_m^{d_m}} \\
& =(k!)^2 \sum_{l,m\leq k}\prod_{i=1}^l \bigg ( \sum_{c_i\geq 1}\sum_{e_i\geq 0}\sum_{x<p_i\leq x^3} \frac{1}{p_i^{\max(c_i,e_i)}}\bigg ) \prod_{i=1}^m \bigg (\sum_{d_i\geq 1}\sum_{x<q_i\leq x^3}\frac{1}{q_i^{d_i}}\bigg ) \\
& =(k!)^2 \sum_{l,m\leq k} \prod_{i=1}^l \bigg ( 2\sum_{x<p_i\leq x^3} \frac{1}{p_i}+O(1)\bigg ) \prod_{i=1}^m \bigg (\sum_{x<q_i\leq x^3}\frac{1}{q_i}+O(1)\bigg )\\
& \ll (k!)^2\sum_{l,m\leq k} \big (2\log 3+C\big )^l \big (\log 3+C\big )^m \ll (Ck)^{2k} .
\end{split}\end{equation} 
The contribution of the case where some exponents are larger than 1 in the right hand side of (\ref{dvd})
is still $\ll (Ck)^{2k}$, by a combination of the previous computation and the argument we used in order to bound $S_2$. \\

Now we analyze the error term, which is
$$ R\ll \frac{1}{\log x}\bigg |\sum_{n\leq x^3}\frac{\Lambda(n)n^{-4/\log x}}{n^{1/2+it}}\bigg |+\frac{k}{\varepsilon}. $$
Hence
\begin{equation}\begin{split}\notag
\frac{1}{T\log T}&\int_T^{2T}|R|^{2k}|\zeta|^2dt \ll \frac{1}{(\log x)^{2k}}\frac{1}{T\log T}\int_T^{2T}\bigg |\sum_{n\leq x^3}\frac{\Lambda(n)n^{-4/\log x}}{n^{1/2+it}}\bigg |^{2k}|\zeta|^2dt+  \frac{k^{2k}}{\varepsilon^{2k}} \\
& \ll \frac{1}{(\log x)^{2k}}\frac{1}{T\log T}\int_T^{2T}\bigg |\sum_{n\leq x^3}\frac{\Lambda(n)n^{-4/\log x}}{n^{1/2+it}}\bigg |^{2k}|\zeta|^2dt+ (\log\log\log T)^{2k}k^{2k}.
\end{split}\end{equation}
We now study the first term, with the aim of proving 
\begin{equation}\label{useful}
\frac{1}{T\log T}\int_T^{2T}\bigg |\sum_{n\leq x^3}\frac{\Lambda(n)n^{-4/\log x}}{n^{1/2+it}}\bigg |^{2k}|\zeta|^2dt \ll (Ck)^{2k}  (\log x)^{2k}. \end{equation}
Using our usual approach we get:
\begin{equation}\begin{split}\label{dvx}
\frac{1}{T\log T}\int_T^{2T}\bigg |\sum_{n\leq x^3}&\frac{\Lambda(n)n^{-4/\log x}}{n^{1/2+it}}\bigg |^{2k}|\zeta|^2dt 
\\& \ll \sum_{\substack{p_1^{r_1},\dots,p_k^{r_k}\leq x^3 \\ q_1^{s_1},\dots, q_k^{s_k}\leq x^3}}\frac{(p_1^{r_1}\cdots p_k^{r_k}, q_1^{s_1}\cdots q_k^{s_k})}{p_1^{r_1}\cdots p_k^{r_k} q_1^{s_1}\cdots q_k^{s_k}}\log p_1 \cdots \log p_k \log q_1 \cdots \log q_k.
\end{split}\end{equation}
Once again we start with the case where all the exponents are equal to 1 and we rewrite the sum in the usual way
$$(k!)^2 \sum_{l,m\leq k}\sum_{\substack{c_1,\dots,c_l\geq 1 \\ c_1+\dots+c_l=k }} \sum_{\substack{e_1,\dots,e_l \geq 0 \\ d_1,\dots,d_m\geq 1 \\ e_1+\dots+e_l+d_1+\dots+d_m=k }} \prod_{i=1}^l \Big ( \sum_{p_i}\frac{(\log p_i)^{e_i+c_i}}{p_i^{\max(e_i,c_i)}} \Big ) \prod_{i=1}^m \Big (\sum_{q_i}\frac{(\log q_i)^{d_i}}{q_i^{d_i}} \Big ).$$
In the case  $\max(e_i,c_i)>1$ (or $d_i>1$) the sum in the first (or second, respectively) parentheses is bounded because of the usual argument. The largest contribution comes from the case $0\leq e_i\leq 1$, $c_i=1$, $d_i=1$, which gives
\begin{equation}\begin{split}\notag
(k!)^2\sum_{l,m\leq k}&\sum_{\substack{0\leq e_1,\dots,e_l \leq 1 \\ e_1+\dots+e_l+m=k }} \prod_{i=1}^l \Big ( \sum_{p_i}\frac{(\log p_i)^{1+e_i}}{p_i} \Big ) \prod_{i=1}^m \Big (\sum_{q_i}\frac{\log q_i}{q_i} \Big ) \\
& \leq (k!)^2\sum_{l,m\leq k} \sum_{\substack{0\leq e_1,\dots,e_l \leq 1 \\ e_1+\dots+e_l+m=k }} \prod_{i=1}^l (3\log x+O(1))^{1+e_i} \prod_{i=1}^m (3\log x+O(1)) \\
& = (k!)^2\sum_{l,m\leq k} \sum_{\substack{0\leq e_1,\dots,e_l \leq 1 \\ e_1+\dots+e_l+m=k }} (3\log x+O(1))^{l+\sum_{i=1}^le_i+m} 
\ll (k!)^2\sum_{l,m\leq k}2^l (4\log x)^{l+k}
\end{split} \end{equation}
and this is $\ll (\log x)^{2k}(Ck)^{2k}$. As before, if some exponents among the $r_i$,$s_j$ in (\ref{dvx}) are larger than 1, then the contribution of this case in (\ref{dvx}) is still $\ll (\log x)^{2k}(Ck)^{2k}$, by a combination of the previous computation and the technique we used to study $S_2$. This proves (\ref{useful}) and as a consequence we get
\begin{equation}\label{difficultbound}
\frac{1}{T\log T}\int_T^{2T} \Bigg (\frac{1}{\log x}\Big (\Big |\sum_{n\leq x^3}\frac{\Lambda(n)n^{-4/\log x}}{n^{1/2+it}} \Big |+\log T \Big ) \Bigg )^{2k}|\zeta|^2dt \ll  (\log\log\log T)^{2k}(Ck)^{2k}.
\end{equation}

What remains to investigate is the contribution of $L(t)$. Following Tsang (\cite{Tsang}, equation (5.21)) we have: 
\begin{equation}\label{L(t)} \Re L(t) \ll L_1(t)+L_2(t) \end{equation}
where denoting with $\rho=\frac{1}{2}+i\gamma$ the non-trivial zeros of $\zeta$ 
$$ L_1(t):=\sum_\rho \Big (\frac{4}{\log x}\Big )^2 \frac{1}{ |\frac{4}{\log x}+i(t-\gamma) |^2}\int_{1/2}^{1/2+4/\log x}\frac{ |u-\frac{1}{2} |}{(u-\frac{1}{2})^2+(t-\gamma)^2} du$$  $$ L_2(t):= \Big (\frac{4}{\log x}\Big )^2 \sum_\rho \frac{1}{|\frac{4}{\log x}+i(t-\gamma)|^2} $$ 
so we need to study the weighted moments of $L_1(t)$ and $L_2(t)$. \\
The latter is not difficult; indeed Selberg proved that (see \cite{Tsang}, equation (5.20))
\begin{equation}\label{Selberg} \sum_\rho \frac{1}{|\frac{4}{\log x}+i(t-\gamma)|^2} \ll \log x \Bigg ( \Big | \sum_{n\leq x^3}\frac{\Lambda(n)n^{-4/\log x}}{n^{1/2+it}} \Big | +\log T \Bigg ) \end{equation}
hence in view of (\ref{difficultbound}) we know that the $2k-$th moment of $L_2(t)$ is $\ll (\log\log\log T)^{2k}(Ck)^{2k}$. 
To deal with $L_1(t)$, we denote $\eta_t:= \min_\rho |t-\gamma|$ and $\log ^+ t:=\max(\log t,0)$. From Tsang's computation (\cite{Tsang}, p.93) we know that:
\begin{multline}\notag L_1(t)\ll \frac{1}{\log x} \Bigg ( \Big | \sum_{n\leq x^3}\frac{\Lambda(n)n^{-4/\log x}}{n^{1/2+it}} \Big |+\log T\Bigg ) \\ +\frac{1}{\log x} \log ^+ \Big(\frac{1}{\eta _t \log x}\Big)  \Bigg ( \Big | \sum_{n\leq x^3}\frac{\Lambda(n)n^{-4/\log x}}{n^{1/2+it}} \Big |+\log T\Bigg )\end{multline} 
and the first term here is not a problem for the same reason as before. As a last step we study the $2k-$th moment of the second term. Applying the Cauchy-Schwarz inequality:
\begin{equation}\begin{split}\label{ultimostep}
\frac{1}{(\log x)^{2k}}&\int_T^{2T} \Big (\log ^+ \frac{1}{\eta _t \log x} \Big )^{2k}\Big ( \Big | \sum_{n\leq x^3} \frac{\Lambda(n)n^{-4/\log x}}{n^{1/2+it}}\Big |+\log T\Big )^{2k}|\zeta|^2dt \\
& \ll \sqrt{T\log T}(\log\log\log T)^{2k}(Ck)^{2k}\sqrt{\int_T^{2T}\Big (\log ^+ \frac{1}{\eta _t \log x} \Big )^{4k}|\zeta|^2dt}.
\end{split}\end{equation}
The proposition follows if we bound the remaining integral. Here the Riemann Hypothesis plays a central role, in the form of a result due to Gonek, which is a consequence of the Landau-Gonek formula (see \cite{Gonek}, Corollary 2 or \cite{Gonek1}, p93, Theorem 2):
\begin{lemma}\label{Gonek}(Gonek)
Assume RH, let T be a large parameter, $\alpha\in\mathbb R$ such that $|\alpha|\leq \frac{\log T}{2\pi}$, then 
$$ \sum_{0<\gamma \leq T}\bigg | \zeta \bigg (\frac{1}{2}+ i\Big (\gamma +\frac{2\pi \alpha}{\log T}\Big ) \bigg )\bigg |^2 = \Big (1-\Big (\frac{\sin \pi\alpha}{\pi \alpha} \Big )^2\Big )\frac{T}{2\pi}(\log T)^2+O\big (T (\log T)^{7/4}\big ).$$ 
\end{lemma}
\noindent For us the uniform upper bound $\ll T(\log T)^2$ for $|\alpha| \leq\frac{\log T}{2\pi \log x}$ will be sufficient.
Using this result we get:
\begin{equation}\begin{split}\notag
\int_T^{2T} \Big (&\log ^+ \frac{1}{\eta _t \log x}\Big )^{4k}  |\zeta|^2dt \\ &\leq \sum_{T-\frac{1}{\log x}\leq \gamma\leq 2T+\frac{1}{\log x}} \int_0^{1/\log x} \Big (\log ^+ \frac{1}{w \log x}\Big )^{4k} \big |\zeta \big (1/2+i(w+\gamma)\big )\big |^2dw \\ & = \sum_{T-\frac{1}{\log x}\leq \gamma\leq 2T+\frac{1}{\log x}} \int_0^1 \Big (\log ^+ \frac{1}{t}\Big )^{4k}\big |\zeta \big (1/2+i(\gamma+\frac{t}{\log x})\big )\big |^2\frac{dt}{\log x} \\ & =\frac{1}{\log x}\int_0^{1} (\log t)^{4k} \sum_{T-\frac{1}{\log x}\leq \gamma\leq 2T+\frac{1}{\log x}}\big |\zeta \big (1/2+i(\gamma+\frac{t}{\log x})\big )\big |^2dt \\ &\ll \frac{T(\log T)^2}{\log x}\int_0^1 (\log t)^{4k}dt \ll T\log T k\varepsilon^{-1}(Ck)^{4k},
\end{split}\end{equation}
since $\int_0^1 (\log t)^{4k}dt=\int_0^\infty e^{-t}t^{4k}dt=\Gamma(4k+1)=(4k)!\ll (4k)^{4k}$. Putting this into (\ref{ultimostep}) one has that also the $2k-$th moment of $L_1(t)$ is bounded by $ (Ck)^{4k}\varepsilon^{-1/2}(\log\log T)^{2k}$. Then the contribution of the zeros is under control, being $$ \frac{1}{T\log T}\int_T^{2T}|\Re L(t)|^{2k}dt\ll (Ck)^{4k}(\log\log\log T)^{2k+1/2} $$ and the proposition follows.

\section{Proof of Proposition \ref{prop2}}
\subsection{Sketch of the proof}
In order to prove Proposition \ref{prop2}, we need to perform a precise asymptotic analysis  for the moments of $\Re P(t)$. First of all, since the polynomial is short ($n\leq x= T^{\varepsilon/k}=T^{o(1/k)}$) one can easily compute its mean and variance by standard applications of Lemma \ref{B-C-HB}. Indeed for any $r,s$ integers one has
\begin{equation}\begin{split}\label{partenza}
\int_T^{2T}&P(t)^r\overline {P(t)}^s|\zeta|^2dt \\ &= T\sum_{\substack{p_1,\dots,p_r\leq x \\ q_1,\dots,q_s\leq x}}\frac{(p_1\cdots p_r,q_1\cdots q_s)}{p_1\cdots p_rq_1\cdots q_s}\bigg (\log \bigg(\frac{T(p_1\cdots p_r,q_1\cdots q_s)^2}{p_1\cdots p_rq_1\cdots q_s}\bigg)+c \bigg )+o(T)
\end{split}\end{equation}
then, since $2\Re P(t)=P(t)+\overline{P(t)}$, the mean of $\Re P(t)$ is
\begin{equation}\begin{split}\notag
\frac{1}{T\log T}\int_T^{2T}\Re P(t)|\zeta|^2dt&=\frac{1}{\log T}\sum_{p\leq x} \frac{\log T-\log p+c}{p}+o\Big(\frac{1}{\log T}\Big) \\ &=  \mathcal L -\frac{\varepsilon}{k} +O\Big(\frac{\log\log T}{\log T}\Big)= \mathcal L+o(1).
\end{split}\end{equation}
Similarly 
\begin{equation}\begin{split}\notag
&\frac{1}{T\log T}\int_T^{2T}(\Re P(t))^2|\zeta|^2dt= \frac{1}{T\log T}\int_T^{2T}\bigg(\frac{1}{4}P(t)^2+\frac{1}{2}P(t)\overline{P(t)}+\frac{1}{4}\overline{P(t)}^2\bigg)|\zeta|^2dt \\
&\;\;= \frac{1}{2\log T}\bigg[\sum_{p_1,p_2\leq x} \frac{\log T-\log (p_1p_2)+c}{p_1p_2} + \sum_{p,q\leq x} \frac{(p,q)}{pq}\bigg (\log \Big(\frac{T(p,q)^2}{pq}\Big)+c\bigg )\bigg] +o\Big(\frac{1}{\log T}\Big) \\
&\;\;=\frac{1}{2\log T}\Big(2\mathcal L^2\log T-4\mathcal L\log x+\mathcal L \log T+O(\log T)\Big)
=\mathcal L^2 + \frac{\mathcal L}{2}-\frac{2\varepsilon\mathcal L}{k}+O(1).
\end{split}\end{equation}
Hence the variance is
\begin{equation}\begin{split}\notag
\frac{1}{T\log T}\int_T^{2T}(\Re P(t)-\mathcal L)^2|\zeta|^2dt\sim \frac{\mathcal L}{2}. 
\end{split}\end{equation}

To prove Proposition \ref{prop2}, we now have to compute the $k-$th moment of $\Re P(t)-\mathcal L$ with respect to $|\zeta|^2dt$, for every $k$ integer.
Here we give a simplified sketch of the proof, leaving the rigorous one for the following section. First of all, since
\begin{equation}\label{approx} \log\bigg(\frac{T(p_1\cdots p_r,q_1\cdots q_s)^2}{p_1\cdots p_rq_1\cdots q_s}\bigg)+c = \log T+\log\bigg(\frac{(p_1\cdots p_r,q_1\cdots q_s)^2}{p_1\cdots p_rq_1\cdots q_s}\bigg)+c\end{equation}
then expanding out the $k-$th power and using (\ref{partenza}) one has
\begin{equation}\begin{split}\label{parto}
\frac{1}{T\log T}\int_T^{2T}&(\Re P(t)-\mathcal L)^k|\zeta|^2dt \\ &=\sum_{j+h=k}\binom{k}{h}(-1)^j\mathcal L^j 2^{-h}\sum_{r+s=h}\binom{h}{r}\sum_{\substack{p_1,\dots,p_r\leq x \\ q_1,\dots,q_s\leq x}}\frac{(p_1\cdots p_r,q_1\cdots q_s)}{p_1\cdots p_rq_1\cdots q_s}+\cdots
\end{split}\end{equation}
where the dots come from the contributions of the second and third terms in (\ref{approx}), which we are going to ignore in the following. Indeed the contribution of the constant $c$ is clearly analogous but smaller than the one coming from $\log T$. Even though the second term in (\ref{approx}) is not negligible compared to the first one, its contribution in the right hand side of (\ref{parto}) can be computed in a similar way to the contribution of the first one, with the important difference that in this case the main term will cancel out. Thus we ignore it as well for now, focusing on the first term.\\
Let's suppose now that the primes $p_1,\dots,p_r$ are distinct and the primes $q_1,\dots,q_s$ are distinct as well. In order to compute explicitly the GCD, we fix an integer $m$, which is smaller than both $r$ and $s$, and we suppose that $m$ repetitions occur among the $p_i$ and the $q_j$. Because of the previous assumptions, it can happen in $\binom{r}{m}\binom{s}{m}m!$ ways (selecting $m$ primes among the $p_i$ and $m$ primes among the $q_j$, then permuting the two blocks multiplying by $m!$), hence
\begin{equation}\begin{split}\notag 
\sum_{\substack{p_1,\dots,p_r \leq x \text{ distinct} \\ q_1,\dots,q_s \leq x \text{ distinct}}} \frac{(p_1\cdots p_r,q_1\cdots q_s)}{p_1\cdots p_rq_1\cdots q_s} = \sum_{m\leq \text{ min}(r,s)}\binom{r}{m}\binom{s}{m}m!\sum_{\substack{p_1,\dots,p_{r+s-m}\leq x \\ \text{distinct}}}\frac{1}{p_1\cdots p_{r+s-m}}.
\end{split}\end{equation}
We now drop the condition in the inner sum that the primes are distinct. As we will show in the following section, all these assumptions about distinct primes do not affect the asymptotic of the moment we are interested in. Indeed the errors coming from all these extra assumptions will all cancel out and give a contribution which is negligible with respect to the main term. With this assumption the previous sum becomes
$$ \sum_{m\leq \text{ min}(r,s)}\frac{1}{m!}\frac{r!}{(r-m)!}\frac{s!}{(s-m)!}\mathcal L^{r+s-m}. $$
Putting this into (\ref{parto}), recalling that $r!/(r-m)!=\partial_X^m[X^r]_{X=1}$, for $k$ even we get:
\begin{equation}\begin{split}\notag
&\frac{1}{T\log T}\int_T^{2T}(\Re P(t)-\mathcal L)^k|\zeta|^2dt \\ 
&\hspace{1cm}= \sum_{j+h=k}\binom{k}{h}(-1)^j\mathcal L^j 2^{-h}\sum_{r+s=h}\binom{h}{r}\sum_{m\leq \text{ min}(r,s)}\frac{1}{m!}\frac{r!}{(r-m)!}\frac{s!}{(s-m)!}\mathcal L^{r+s-m}+\cdots \\
&\hspace{1cm}=\sum_{m\leq \frac{k}{2}}\frac{\mathcal L^{k-m}}{m!}\Big [\frac{k!}{(k-2m)!}\Big(\frac{X+Y}{2}-1 \Big)^{k-2m} 2^{-2m} \Big ]_{X=Y=1} +\cdots\\ 
&\hspace{1cm}=\sum_{m\leq \frac{k}{2}}\frac{\mathcal L^{k-m}}{2^{2m}m!}\frac{k!}{(k-2m)!} \mathbbm{1}_{2m=k}+\cdots = 
\frac{k!}{2^k(k/2)!}\mathcal L^{k/2} +\cdots= \Big (\frac{\mathcal L}{2}\Big )^{k/2}(k-1)!!+\cdots
\end{split}\end{equation}
since $k!=2^{k/2}(k/2)!(k-1)!!$ for any even $k$. Otherwise if $k$ is odd, then the main term vanishes, being $m\leq (k-1)/2$. \\

We now highlight the main difference from the classical case \cite{RS}. There one easily sees that $\int P(t)^r\overline{P(t)}^sdt$ is non negligible only if $r$ equals $s$. Therefore just the diagonal term $r=s=k/2$ contributes to the main term of the $k-$th moment of $\Re P(t)$. On the other hand this is no longer true in the weighted case, since all the integrals $\int P(t)^r\overline{P(t)}^s|\zeta|^2dt$ give a contribution of order $T\log T\mathcal L^{r+s}$. The main point is that in the classical case the mean of $\Re P(t)$ is 0, while with respect to the weighted measure $|\zeta|^2dt$ the mean is $\sim\mathcal L$. Thus, even though in the weighted case the size of the $k-$th moment of $\Re P(t)$ is $\mathcal L^k$, the $k-$th moment of $\Re P(t)-\mathcal L$ has order $\mathcal L^{k/2}$. Showing this cancellation from $k$ to $k/2$ is the bulk of the proof.
\subsection{Proof of Proposition \ref{prop2}}

We now prove the result, following the line of the previous computation. Expanding out the $k-$th power and using $2\Re P(t)=P(t)+\overline{P(t)}$, one finds
\begin{equation}\begin{split}\label{start}
\int_T^{2T}(\Re P(t)&-\mathcal L)^k|\zeta|^2dt =\sum_{j+h=k}\binom{k}{h}(-1)^j\mathcal L^j 2^{-h}\sum_{r+s=h}\binom{h}{r}\int_T^{2T}P(t)^r\overline{P(t)}^s|\zeta|^2dt
\end{split}\end{equation}
and the inner integral equals
$$T\sum_{\substack{p_1,\dots,p_r\leq x \\ q_1,\dots,q_s\leq x}}\frac{(p_1\cdots p_r,q_1\cdots q_s)}{p_1\cdots p_rq_1\cdots q_s}\bigg (\log\bigg(\frac{T(p_1\cdots p_r,q_1\cdots q_s)^2}{p_1\cdots p_rq_1\cdots q_s}\bigg)+c\bigg )+o(T)$$
in view of (\ref{partenza}). Since $\log t=\partial_w [t^w]_{w=0}$, one gets
\begin{equation}\begin{split}\label{rs}
\int_T^{2T}P(t)^r\overline{P(t)}^s|\zeta|^2dt=T(\log T+c) f_x(0)+T\partial_w[f_x(w)]_{w=0}+o(T)
\end{split}\end{equation}
where 
\begin{equation}\begin{split}\notag
f_x(w)=\sum_{\substack{p_1,\dots,p_r\leq x \\ q_1,\dots,q_s\leq x}}\frac{(p_1\cdots p_r,q_1\cdots q_s)^{2w+1}}{(p_1\cdots p_rq_1\cdots q_s)^{w+1}}
\end{split}\end{equation}
In order to be able to compute explicitly the GCD, we put in evidence the possible repetitions among the primes, re-writing the $p_i$ and the $q_i$ as follows. First we put in evidence the repetitions among the primes $p_i$, writing
$$ p_1,\dots,p_r \longrightarrow p_1,\dots,p_{r-v_1},p{'}_{1} ^{\alpha_1},\dots,p{'}_{u_1} ^{\alpha_{u_1}} $$
where $p_1,\dots,p_{r-v_1},p{'}_{1} ,\dots,p{'}_{u_1}$ are all distinct, $\alpha_1+\dots+\alpha_{u_1}=v_1$, $\alpha_i\geq 2$ for every $i$. With this change of variable we need a normalization $ \frac{r!}{(r-v_1)!}c_{\underline{\alpha}} $, where $c_{\underline\alpha}$ is a positive coefficient smaller than 1, which does not depend on $r$ but just on the configuration $\alpha_1,\dots,\alpha_{u_1}$. Notice that if $v_1=0$, then $c_{\underline \alpha}=1$. 
Now we highlight the multiplicities of the primes $q_j$ and we put in evidence those ones that already appear among the $p'_i$. Then we write
$$ q_1,\dots,q_s \longrightarrow q_1,\dots,q_{s-v_2-a_2},p'_1,\dots,p'_{a_2},q{'}_{1} ^{\beta_1},\dots,q{'}_{u_2} ^{\beta_{u_2}},p{'}_1^{\gamma_1},\dots,p{'}_{u_1}^{\gamma_{u_1}} $$
with  $q_i$ distinct, $q'_i$ distinct, $q'_i\not =p{'}_j$ for every $i,j$, $q_i\not = q'_j,p'_j$ for every $i,j$ and $\beta_1+\dots+\beta_{u_2}+\gamma_1+\dots+\gamma_{u_1}+a_2=v_2$, $\beta_i\geq 2$, $\gamma_i\not =1$ for every $i$. Also in this case the change of variable brings into play a normalization $\binom{s-v_2}{a_2}\binom{u_1}{a_2}a_2!\frac{s!}{(s-v_2)!}c_{\underline{\beta},\underline{\gamma}}$, where once again $c_{\underline{\beta},\underline{\gamma}}$ only depends on the configuration $\beta_1,\dots,\beta_{u_2},\gamma_1,\dots,\gamma_{u_1}$ and it is equal to $1$ when $u_2=0$ and $\gamma_i=0$ for every $i$. 
The normalization coefficient comes from standard combinatorics as follows. We make the multiplicity of any $q_i$ explicit, putting in evidence the $s-v_2$ ones which appear once. This can be done in $s!/(s-v_2)!$ ways times a coefficient described above, which does not depend on $s$. Moreover, in order to highlight the coincidences between the $q_i$ and the $p'_i$ (say we have $a_2$ coincidences), we select $a_2$ primes among $q_1,\dots,q_{s-v_2}$ ($\binom{s-v_2}{a_2}$ ways) and $a_2$ primes among $p'_1,\dots p'_{u_1}$ ($\binom{u_1}{a_2}$ ways) and then we permute the two blocks multiplying by $a_2!$.
Then we have
\begin{equation}\begin{split}\notag
f_x(w)=&\sum_{\substack{v_1,u_1\leq r \\ v_2,u_2\leq s}}\sum_{a_2\leq s-v_2}\binom{s-v_2}{a_2}\binom{u_1}{a_2}a_2!\sum_{\substack{\alpha_1+\dots+\alpha_{u_1}=v_1 \\ \alpha_i\geq 2\;\;\forall i}}c_{\underline\alpha}\sum_{\substack{\beta_1+\dots+\beta_{u_2}+\gamma_1+\dots+\gamma_{u_1}=v_2 \\ \beta_i\geq 2,\; \gamma_i \not = 1\;\;\forall i}}c_{\underline\beta ,\underline\gamma}
\\ &\frac{r!}{(r-v_1)!}\frac{s!}{(s-v_2)!}\sum_{\substack{p'_i,q'_j \\ \text{distinct}}}\frac{(p{'}_1^{\alpha_1}\cdots p{'}_{u_1}^{\alpha_{u_1}},p'_1\cdots p'_{a_2}p{'}_1^{\gamma_1}\cdots p{'}_{u_1}^{\gamma_{u_1}})^{2w+1}}{(p{'}_1^{\alpha_1}\cdots p{'}_{u_1}^{\alpha_{u_1}}p'_1\cdots p'_{a_2}q{'}_1^{\beta_1}\cdots q{'}_{u_2}^{\beta_{u_2}} p{'}_1^{\gamma_1}\cdots p{'}_{u_1}^{\gamma_{u_1}})^{w+1}} \\
& \sum_{\substack{p_1,\dots,p_{r-v_1}\leq x \\ \text{distinct and }\not = p'_i \\ q_1,\dots,q_{s-v_2-a_2}\leq x \\ \text{distinct and }\not = p'_i,q'_j}}\frac{[(p_1\cdots p_{r-v_1},q_1\cdots q_{s-v_2-a_2})(p_1\cdots p_{r-v_1},q{'}_1^{\beta_1}\cdots q{'}_{u_2}^{\beta_{u_2}})]^{2w+1}}{(p_1\cdots p_{r-v_1}q_1\cdots q_{s-v_2-a_2})^{w+1}}.
\end{split}\end{equation}
For the sake of brevity let's denote $\underline p'$ and $\underline q'$ the product of $p'_i$ and $q'_i$ respectively with their exponents (for instance $\underline p'^{\underline\alpha}:=p{'}_1^{\alpha_1}\cdots p{'}_{u_1}^{\alpha_{u_1}}$). To be able to compute the GCD between $\underline p$ and $\underline q'^{\underline\beta}$ in the inner sum, we now put in evidence the repetitions among the $p_i$ and the $q_j'$. Let's say we have $a_1$ primes among the $p_i$ which coincide with some $q_j'$. Then, denoting $r':=r-v_1-a_1$ and $s':=s-v_2-a_2$, we get
\begin{equation}\begin{split}\notag
f_x(w)=&\sum_{\substack{v_1,u_1\leq r \\ v_2,u_2\leq s}}\sum_{\substack{a_1\leq r-v_1 \\ a_2\leq s-v_2}}\sum_{\underline{\alpha},\underline{\beta},\underline{\gamma}}c(\underline\alpha,\underline{\beta},\underline{\gamma},\underline a) \frac{r!}{(r-v_1-a_1)!}\frac{s!}{(s-v_2-a_2)!}
\\ &\sum_{\substack{p'_i,q'_j \\ \text{distinct}}}\frac{(\underline p{'}^{\underline{\alpha}},\underline p' \underline p{'}^{\underline \gamma})^{2w+1}(\underline q')^w}{(\underline p{'}^{\underline \alpha}\underline p'\underline q{'}^{\underline\beta}\underline p{'}^{\underline \gamma})^{w+1}} \sum_{\substack{p_1,\dots,p_{r'} \text{ distinct and }\not = p'_i,q'_j \\ q_1,\dots,q_{s'} \text{ distinct and }\not = p'_i,q'_j}}\frac{(p_1\cdots p_{r'},q_1\cdots q_{s'})^{2w+1}}{(p_1\cdots p_{r'}q_1\cdots q_{s'})^{w+1}}
\end{split}\end{equation}
where $c(\underline\alpha,\underline{\beta},\underline{\gamma},\underline a)$ is a bounded coefficient which does not depend on $r$ and $s$ and it is equal to 1 when $u_i=v_i=a_i=0$ for $i=1,2$.
Note that the sum over $p'_i$ and $q'_j$ is bounded when $w$ is close to $0$, since both $\beta_i$ and max$(\alpha_i,\gamma_i+1)$ are $\geq 2$. Lastly we want to put in evidence the repetitions among the $p_i$ and the $q_j$, in order compute explicitly the last greatest common divisor $(p_1\cdots p_{r'},q_1\cdots q_{s'})$ in the inner sum. If $m$ repetitions occur, for any $m\leq \min(r',s')$, we finally have $r'+s'-m$ distinct primes and the coefficient of normalization is $\binom{r'}{m}\binom{s'}{m}m!$. Therefore
\begin{equation}\begin{split}\label{inutile}
f_x(w)&=\sum_{\substack{v_1,u_1\leq r \\ v_2,u_2\leq s}}\sum_{\substack{a_1\leq r-v_1 \\ a_2\leq s-v_2}}\sum_{\underline{\alpha},\underline{\beta},\underline{\gamma}}c(\underline\alpha,\underline{\beta},\underline{\gamma},\underline a) \sum_{\substack{p'_i,q'_j \\ \text{distinct}}}\frac{(\underline p{'}^{\underline{\alpha}},\underline p' \underline p{'}^{\underline \gamma})^{2w+1}(\underline q')^w}{(\underline p{'}^{\underline \alpha}\underline p'\underline q{'}^{\underline\beta}\underline p{'}^{\underline \gamma})^{w+1}} 
\\ &\sum_{m\leq \min(r',s')}\frac{r!}{(r'-m)!}\frac{s!}{(s'-m)!}\frac{1}{m!}\sum_{\substack{p_1,\dots,p_{r'+s'-2m} \\ q_1,\dots,q_m \\ \text{ distinct and }\not = p'_i,q'_j }}\frac{1}{(p_1\cdots p_{r'+s'-2m})^{w+1}q_1\cdots q_{m}}.
\end{split}\end{equation}
After computing the GCD, we now remove the extra conditions in the inner sum, which force the primes $p_i$ and $q_j$ to be all distinct and $\not = p'_i,q'_j $. We get rid of the condition that forces the primes to be all distinct by using basic combinatorics and we remove the last condition $p_1,\dots,p_{r'+s'-2m},q_1,\dots,q_m\not = p'_i,q'_j $, splitting the inner sums as 
$$ \sum_{\substack{p\leq x \\ p\not =p'_i,q'_{j}}}\frac{1}{p^s}= \sum_{p\leq x}\frac{1}{p^s}-\sum_{i=1}^{u_1}\frac{1}{p{'}^{s}}-\sum_{i=1}^{u_2}\frac{1}{q{'}^{s}}$$
and expanding out the powers by Newton's binomial formula. Hence we have (denote $h':=r'+s'$)

\begin{equation}\begin{split}\notag
f_x(w)=&
\\&\sum_{\substack{v_1,u_1\leq r \\ v_2,u_2\leq s}}\sum_{\substack{a_1\leq r-v_1 \\ a_2\leq s-v_2}}\sum_{\underline{\alpha},\underline{\beta},\underline{\gamma}} \sum_{\substack{ p'_i, q'_j \\ \text{distinct}}}\frac{(\underline p{'}^{\underline{\alpha}},\underline p' \underline p{'}^{\underline \gamma})^{2w+1}(\underline q')^w}{(\underline p{'}^{\underline \alpha}\underline p'\underline q{'}^{\underline\beta}\underline p{'}^{\underline \gamma})^{w+1}} 
\sum_{m\leq \min(r',s')}\frac{r!}{(r'-m)!}\frac{s!}{(s'-m)!}\\&\sum_{\substack{t_1 \leq h'-2m \\ t_2\leq m \\ t_3 \leq t_1+t_2}}\sum_{\substack{\mathcal {P}\in Part \\ \mathcal{P} =\{R_1,\dots,R_{t_3}\} \\ \text{with } r_i:=\sum_{j\in R_i}a_j \geq 2}} \prod_{i=1}^{t_3}(\#R_i-1)!(-1)^{\#R_i-1}\Big ( \sum_{p\not = \underline p',\underline q'}\frac{1}{p^{r_i}} \Big )c(\underline\alpha,\underline{\beta},\underline{\gamma},\underline a,\underline t) \\ &
\sum_{l_1\leq h'-2m-t_1}\frac{(h'-2m)!}{l_1!(h'-2m-t_1-l_1)!}\Big ( -\sum_{i=1}^{u_1}\frac{1}{p{'}^{1+w}}-\sum_{i=1}^{u_2}\frac{1}{q{'}^{1+w}}\Big )^{l_1}  \\&\Big ( \sum_{p\leq x }\frac{1}{p^{1+w}} \Big )^{h'-2m-t_1-l_1}
\sum_{l_2\leq m-t_2}\frac{1}{l_2!(m-t_2-l_2)!}\Big ( -\sum_{i=1}^{u_1}\frac{1}{p{'}}-\sum_{i=1}^{u_2}\frac{1}{q{'}}\Big )^{l_2}\mathcal L^{m-t_2-l_2}
\end{split}\end{equation}
where $c(\underline\alpha,\underline{\beta},\underline{\gamma},\underline a,\underline t)$ is a bounded coefficient not depending on $r,s,m$, which is equal to 1 if the parameters $v_i,u_i,t_i$ are all equal to 0 and $Part$ denotes the set of partitions of the set of the exponents of primes appearing in the inner sum in (\ref{inutile}).\\ 

We are ready to plug the formula we got for $f_x(w)$ into the formula for the $k-$th moment of $\Re P(t)-\mathcal L$. Putting (\ref{start}) and (\ref{rs}) together one has 
\begin{equation}\begin{split}\label{split}
\int_T^{2T}(\Re& P(t)-\mathcal L)^k|\zeta|^2dt \\ 
&=(T(\log T+c))\Big [\sum_{j+h=k}\binom{k}{h}(-1)^j\mathcal L^j 2^{-h}\sum_{r+s=h}\binom{h}{r} f_x(w)\Big ]_{w=0} \\&\quad +\partial_w \Big [T\sum_{j+h=k}\binom{k}{h}(-1)^j\mathcal L^j 2^{-h}\sum_{r+s=h}\binom{h}{r}f_x(w)\Big ]_{w=0}+o(T)
\end{split}\end{equation}
Now we exchange the order of summation, bringing the sum over $j,h$ inside in order to appreciate the cancellation. By the explicit expression we got for $f_x(w)$ we have 
\begin{equation}\begin{split}\label{casa}
\sum_{j+h=k}\binom{k}{h}&(-1)^j\mathcal L^j 2^{-h}\sum_{r+s=h}\binom{h}{r} f_x(w) 
\\ =&
\sum_{m\leq \frac{k}{2}}\quad \sum_{\substack{\underline v,\underline u,\underline a , \underline \alpha,\underline\beta,\underline \gamma ,\underline t ,\mathcal P ,\underline l\\  p'_i, q'_j\text{ distinct}  }} 
F_{\underline v,\underline u,\underline a , \underline \alpha,\underline\beta,\underline \gamma ,\underline t ,\mathcal P ,\underline l}(p'_i,q'_j;w)
\frac{1}{(m-t_2-l_2)!}\mathcal L^{m-t_2-l_2}
\\& \sum_{j+h=k}\binom{k}{h}(-1)^j\mathcal L^{j}2^{-h}\frac{(h'-2m)!}{(h'-2m-t_1-l_1)!}\Big ( \sum_{p\leq x }\frac{1}{p^{1+w}} \Big )^{h'-2m-t_1-l_1}
\\&  \sum_{r+s=h}\binom{h}{r}\frac{r!}{(r-v_1-a_1-m)!}\frac{s!}{(s-v_2-a_2-m)!}
\end{split}\end{equation}
where we denote $k':=k-v_1-v_2-a_1-a_2$ and
\begin{equation}\begin{split}\notag
F_{\underline v,\underline u,\underline a , \underline \alpha,\underline\beta,\underline \gamma ,\underline t ,\mathcal P ,\underline l}&(p'_i,q'_j;w) := \frac{(\underline p{'}^{\underline{\alpha}},\underline p' \underline p{'}^{\underline \gamma})^{2w+1}(\underline q')^w}{(\underline p{'}^{\underline \alpha}\underline p'\underline q{'}^{\underline\beta}\underline p{'}^{\underline \gamma})^{w+1}} 
c(\underline\alpha,\underline{\beta},\underline{\gamma},\underline a,\underline t,\underline l)\Big ( \sum_{p\not = p'_i,q'_j}\frac{1}{p^{r_i}} \Big ) \\ &
\prod_{i=1}^{t_3}(\#R_i-1)!(-1)^{\#R_i-1}
\Big ( -\sum_{i=1}^{u_1}\frac{1}{p{'}^{1+w}}-\sum_{i=1}^{u_2}\frac{1}{q{'}^{1+w}}\Big )^{l_1}\Big ( -\sum_{i=1}^{u_1}\frac{1}{p{'}}-\sum_{i=1}^{u_2}\frac{1}{q{'}}\Big )^{l_2} .
\end{split}\end{equation}
Note that the function $F_{\underline v,\underline u,\underline a , \underline \alpha,\underline\beta,\underline \gamma ,\underline t ,\mathcal P ,\underline l}(p'_i,q'_j;w)$ makes the sum over $p_i,q_j$ in (\ref{casa}) converge. Moreover notice that in the case trivial case $v_i=u_i=a_i=t_i=l_i=0$ for every $i$ then we have that $F_{\underline v,\underline u,\underline a , \underline \alpha,\underline\beta,\underline \gamma ,\underline t ,\mathcal P ,\underline l}(p'_i,q'_j;w)=1$. Now we recall that the three quotients involving $r!$, $s!$ and $h'!$ can be expressed in terms of derivatives (for instance $r!/(r-v_1-a_1-m)! =\partial_X^{v_1+a_1+m}[X^r]_{X=1}$) then (\ref{casa}) becomes
\begin{equation}\begin{split}\notag
&\sum_{m\leq \frac{k}{2}}\quad \sum_{\substack{\underline v,\underline u,\underline a , \underline \alpha,\underline\beta,\underline \gamma ,\underline t ,\mathcal P ,\underline l\\  p'_i, q'_j\text{ distinct}  }} 
F_{\underline v,\underline u,\underline a , \underline \alpha,\underline\beta,\underline \gamma ,\underline t ,\mathcal P ,\underline l}(p'_i,q'_j;w)
\frac{1}{(m-t_2-l_2)!}\mathcal L^{m-t_2-l_2}
\\&
\partial_X^{v_1+a_1+m}\partial_Y^{v_2+a_2+m}\partial_Z^{t_1+l_1}
\Big [ \Big ( \sum_{p\leq x }\frac{1}{p^{1+w}} \Big )^{-v_1-v_2-a_1-a_2-2m-t_1-l_1}Z^{-v_1-v_2-a_1-a_2-2m}\\&
\sum_{j+h=k}\binom{k}{h}(-1)^j\mathcal L^{j}2^{-h}Z^{h}
\Big ( \sum_{p\leq x }\frac{1}{p^{1+w}} \Big )^{h}
(X+Y)^h \Big ]_{X=Y=Z=1}.
\end{split}\end{equation}
Carrying out the computation straightforwardly, denoting $y=v_1+v_2+a_1+a_2+t_1+l_1$, it yields
\begin{equation}\begin{split}\label{cancellation}
\sum_{j+h=k}\binom{k}{h}(-1)^j\mathcal L^j& 2^{-h}\sum_{r+s=h}\binom{h}{r} f_x(w) \\=&
\sum_{m\leq \frac{k}{2}}\quad \sum_{\substack{\underline v,\underline u,\underline a , \underline \alpha,\underline\beta,\underline \gamma ,\underline t ,\mathcal P ,\underline l\\ p'_i,q'_j\text{ distinct}  }} 
F_{\underline v,\underline u,\underline a , \underline \alpha,\underline\beta,\underline \gamma ,\underline t ,\mathcal P ,\underline l}(p'_i, q'_j;w)
\frac{\mathcal L^{m-t_2-l_2}}{(m-t_2-l_2)!}
\\ & 
\frac{k(k-1)\cdots (k-y-2m+1)}{2^{y-t_1-l_1+2m}}
\bigg (\Big ( \sum_{p\leq x }\frac{1}{p^{1+w}} \Big )-\mathcal L\bigg ) ^{k-y-2m} .
\end{split}\end{equation}
Now, recalling (\ref{split}), we have to study the right hand side of (\ref{cancellation}) and its derivative at $w=0$. As we will see soon, only the former contributes to the main term of the $k-$th moment we are considering.  
\\ 

By definition of $\mathcal L:=\sum_{p\leq x}\frac{1}{p}$, if $w=0$ then the expression in the parentheses on the right hand side of (\ref{cancellation}) vanishes. This forces its exponent to be zero, otherwise all the contribution vanishes. Hence we get
\begin{equation}\begin{split}\label{primo}
\Big [\sum_{j+h=k}\binom{k}{h}(-1)^j\mathcal L^j 2^{-h}&\sum_{r+s=h}\binom{h}{r} f_x(w)\Big ]_{w=0} 
\\ =&
\sum_{m\leq \frac{k}{2}}\quad \sum_{\substack{\underline v,\underline u,\underline a , \underline \alpha,\underline\beta,\underline \gamma ,\underline t ,\mathcal P ,\underline l\\ p'_i,q'_j\text{ distinct}  }} 
F_{\underline v,\underline u,\underline a , \underline \alpha,\underline\beta,\underline \gamma ,\underline t ,\mathcal P ,\underline l}(p'_i,q'_j;w)
\frac{\mathcal L^{m-t_2-l_2}}{(m-t_2-l_2)!}
\\ & \hspace{4.8cm}
\frac{k(k-1)\cdots (k-y-2m+1)}{2^{y-t_1-l_1+2m}}
\mathbbm {1}_{2m=k-y}
\end{split}\end{equation}
The main term is given by the largest $m$ possible, i.e. $m=\frac{k}{2}$ if $k$ is even. Since $2m=k-y$, then $y=0$ hence all the parameters that individuate the configuration vanish. Therefore
\begin{equation}\begin{split}\label{I}
\Big [\sum_{j+h=k}&\binom{k}{h}(-1)^j\mathcal L^j 2^{-h}\sum_{r+s=h}\binom{h}{r} f_x(w)\Big ]_{w=0} = \frac{k!}{2^k(k/2)!}\mathcal L^{k/2}+O_k(\mathcal L^{k/2-1})
\end{split}\end{equation}
which matches with the $k-$th moment of a Gaussian by basic properties of the double factorial, since $k!=2^{k/2}(k/2)!(k-1)!!$ for any even $k$. Note that the error term in (\ref{I}) is given by the term $m=k/2-1$ hence it is $O_k(\mathcal L^{k/2-1})$. Of course if $k$ is odd one can immediately see that the right hand side of (\ref{primo}) is $O_k(\mathcal L^{(k-1)/2})$.\\

Let's now analyze the derivative
\begin{equation}\begin{split}\label{secondo}
\partial_w \Big [\sum_{j+h=k}&\binom{k}{h}(-1)^j\mathcal L^j 2^{-h}\sum_{r+s=h}\binom{h}{r}f_x(w)\Big ]_{w=0} 
\\=&
\partial_w \Big [
\sum_{m\leq \frac{k}{2}}\quad \sum_{\substack{\underline v,\underline u,\underline a , \underline \alpha,\underline\beta,\underline \gamma ,\underline t ,\mathcal P ,\underline l\\  p'_i,q'_j\text{ distinct}  }} 
F_{\underline v,\underline u,\underline a , \underline \alpha,\underline\beta,\underline \gamma ,\underline t ,\mathcal P ,\underline l}(p'_i,q'_j;w)
\frac{1}{(m-t_2-l_2)!}\mathcal L^{m-t_2-l_2}
\\ &
 \frac{k(k-1)\cdots (k-y-2m+1)}{2^{y-t_1-l_1+2m}}
\bigg (\Big ( \sum_{p\leq x }\frac{1}{p^{1+w}} \Big )-\mathcal L \bigg ) ^{k-y-2m} 
\Big ]_{w=0} .
\end{split}\end{equation}
Recall that this term will be multiplied by a factor $T$ in (\ref{split}), while the other one by $T\log T$. When we compute the derivative using Leibniz's rule, the term where the derivative of $F$ appears is trivially $O_k(\mathcal L^k/\log T)$, which is negligible. Indeed the sum over $p'_i, q'_j$ is still bounded because the exponents of the variables are larger that $2$ and computing derivatives just $\log p_i$ or $\log q_j$ come out. 
We finally have to deal with the derivative of the inner term. Since
$$ \partial_w \Big [ \sum_{p\leq x }\frac{1}{p^{1+w}} \Big ]_{w=0}=
-\sum_{p\leq x}\frac{\log p}{p}\ll \log x=\frac{\varepsilon}{k}\log T$$
we get that the contribution coming from derivative of $p^{-1-w}$ in (\ref{secondo}) is
\begin{equation}\begin{split}\notag
\ll\sum_{m\leq \frac{k}{2}}& \sum_{\substack{\underline v,\underline u,\underline a \\ \underline \alpha,\underline\beta,\underline \gamma ,\underline t ,\mathcal P ,\underline l\\ p'_i,q'_j\text{ distinct}  }} 
F_{\underline v,\underline u,\underline a , \underline \alpha,\underline\beta,\underline \gamma ,\underline t ,\mathcal P ,\underline l}(p'_i,q'_j;0)
\frac{1}{(m-t_2-l_2)!}\mathcal L^{m-t_2-l_2}
\\ &
 \frac{k(k-1)\cdots (k-y-2m+1)}{2^{y-t_1-l_1+2m}}
\partial_w\Bigg [\bigg ( \sum_{p\leq x }\frac{1}{p^{1+w}}-\mathcal L \bigg ) ^{k-y-2m} 
\Bigg ]_{w=0} 
\end{split}\end{equation}
which is $ O_k\big(\varepsilon\log T \mathcal L^{(k-1)/2}\big) $ by the same argument as before. Hence 
\begin{equation}\begin{split}\label{II}
\partial_w \Big [\sum_{j+h=k}&\binom{k}{h}(-1)^j\mathcal L^j 2^{-h}\sum_{r+s=h}\binom{h}{r}f_x(w)\Big ]_{w=0} =O_k\big (\varepsilon\log T \mathcal L^{(k-1)/2} \big )
\end{split}\end{equation}
Putting both (\ref{I}) and (\ref{II}) into (\ref{split}) the proof is complete. \\ 

\textbf{Acknowledgments}.  I am grateful to my advisor Sandro Bettin for his support and encouragement. I would also like to thank Maksym Radziwill for suggesting me this problem and for some helpful conversations and Adam Harper for pointing me out the connection between Theorem \ref{thm1} and Girsanov's theorem.


{\small

\end{document}